\documentclass[10pt,draft,numreferences]{kluwer}
\usepackage[british]{babel}

\newcommand*{\RR}{\mathbb{R}}
\newcommand*{\CC}{\mathbb{C}}
\newcommand*{\NN}{\mathbb{N}}
\newcommand*{\fM}{\mathfrak{M}}
\newcommand*{\fMe}{\fM_1}
\newcommand{\dd}{\:{\mathstrut\mathrm{d}}}
\newcommand*{\motimes}{\otimes}
\def\supp{\mathop{\:\mathrm{supp}\:}}

\newtheorem{proposition}{Proposition}
\newtheorem{theorem}[proposition]{Theorem}
\newtheorem{corollary}[proposition]{Corollary}
\newtheorem{lemma}[proposition]{Lemma}

\newdisplay{example}[proposition]{Example}

\newproof{definition}{Definition}
\newproof{remark*}{Remark}

\newenvironment{equation*}{$$}{$$}
\def\text#1{\mbox{\scriptsize #1}}
\def\eqref#1{(\ref{#1})}

\def\qedsymbol{\rule{1ex}{1ex}}
\def\qed{\rule{0pt}{0pt}\hfill\qedsymbol}
\makeatletter
\def\@copyrightowner{}

\makeatother
\begin{document}
\begin{opening}
\title{Transfinite diameter, Chebyshev constant and
energy on locally compact spaces}

\author{B\'{a}lint Farkas\thanks{This work was started during the 3$^{\text{rd}}$
Summerschool on Potential Theory, 2004, hosted by the College of
Kecskem\'{e}t, Faculty of Mechanical Engineering and Automation
(GAMF). Both authors would like to express their gratitude for the
hospitality and the support received during their stay in
Kecskem\'{e}t.} \email{farkas@mathematik.tu-darmstadt.de}}
\institute{Technische Universit\"{a}t Darmstadt, Fachbereich
Mathematik\newline Department of Applied Analysis\newline
Schlo\ss{}gartenstra\ss{}e 7, D-64289, Darmstadt, Germany}
\author{B\'{e}la Nagy\thanks{The second named author was supported by the Hungarian Scientific Research Fund; OTKA 49448}\email{nbela@sol.cc.u-szeged.hu}}
\institute{Bolyai Institute, University of Szeged\newline Aradi
v\'{e}rtan\'{u}k tere 1\newline H-6720, Szeged, Hungary}

\dedication{Dedicated to the memory of Professor Gustave Choquet (1
March 1915 - 14 November 2006)}

\date{}

\classification{2000 Math. Subj. Class.}{31C15; 28A12, 54D45}


\begin{abstract}
We study the relationship between transfinite diameter, Chebyshev
constant and Wiener energy in the abstract linear potential analytic
setting pioneered by Choquet, Fuglede and Ohtsuka. It turns out
that, whenever the potential theoretic kernel has the maximum
principle, then all these quantities are equal for all compact sets.
For continuous kernels even the converse statement is true: if the
Chebyshev constant of any compact set coincides with its transfinite
diameter, the kernel must satisfy the maximum principle. An
abundance of examples is provided to show the sharpness of the
results.
\end{abstract}

\end{opening}


\section{Introduction}
\label{sec:Intro}

The idea behind  abstract (linear) potential theory, as developed by
Choquet \cite{choquet:1958}, Fuglede \cite{fuglede:1960} and Ohtsuka
\cite{ohtsuka:1961}, is to replace the Euclidian space $\RR^d$ by
some locally compact space $X$ and the well-known Newtonian kernel
by some other kernel function $k:X\times X\to \RR\cup\{+\infty\}$,
and to look at which ``potential theoretic'' assertions remain true
in this generality (see the monograph of Landkof
\cite{landkof:1992}). This approach facilitates general
understanding of certain potential theoretic phenomena and allows
also the exploration of fundamental principles like Frostman's
maximum principle.

Although there is a vast work done considering energy integrals and
different notions of energies, the familiar notions of transfinite
diameter and Chebyshev constants in this abstract setting are
sporadically found, sometimes indeed inaccessible, in the
literature, see Choquet \cite{choquet:1958} or Ohtsuka
\cite{ohtsuka:1967}. In \cite{choquet:1958} Choquet defines
\emph{transfinite diameter} and proves its equality with the
\emph{Wiener energy} in a rather general situation, which of course
covers the classical case of the logarithmic kernel on $\CC$. We
give a slightly different definition for the transfinite diameter
that, for infinite sets, turns out to be equivalent with the one of
Choquet. The primary aim of this note is to revisit the above
mentioned notions and related results and also to partly complement
the theory.

We already remark here that Zaharjuta's generalisation of
transfinite diameter and Chebyshev constant to $\CC^n$ is completely
different in nature, see \cite{zaharjuta:1975}, whereas some
elementary parts of weighted potential theory (see, e.g., Mhaskar,
Saff \cite{mhaskar/saff:1992} and Saff, Totik
\cite{saff/totik:1997}) could fit in this framework.

\medskip

The power of the abstract potential analytic tools is well
illustrated by the notion of the average distance number  from
metric analysis, see Gross \cite{gross:1964}, Stadje
\cite{stadje:1981}. The surprising phenomenon noticed by Gross is
the following: If $(X,d)$ is a compact connected metric space, there
always exists a unique number $r(X)$ (called the \emph{average
distance number} or the \emph{rendezvous number} of $X$), with the
property that for any finite point system $x_1,\dots, x_n\in X$
there is another point $x\in X$ with average distance
\begin{equation*}
\frac{1}{n}\sum_{j=1}^n d(x_j,x)=r(X).
\end{equation*}
Stadje generalised this to arbitrary continuous, symmetric functions
replacing $d$.  Actually, it turned out, see the series of papers
\cite{farkas/revesz:2004c,farkas/revesz:2004a,farkas/revesz:2004b}
and the references therein, that many of the known results
concerning average distance numbers (existence, uniqueness, various
generalisations, calculation techniques etc.), can be proved in a
unified way using the works of Fuglede and Ohtsuka. We mention for
example that Frostman's Equilibrium Theorem is to be accounted for
the existence for certain invariant measures (see Section
\ref{sec:sum} below). In these investigations the two variable
versions of Chebyshev constants and energies and even their minimax
duals had been needed, and were also partly available due to the
works of Fuglede \cite{fuglede:1965} and Ohtsuka
\cite{ohtsuka:1965,ohtsuka:1967}, see also
\cite{farkas/revesz:2004c}.

Another occurrence of abstract Chebyshev constants is in the study
of polarisation constants of normed spaces, see Anagnostopoulos,
R\'e\-v\'esz \cite{anagnostopoulos/revesz:2006} and R\'ev\'esz,
Sarantopoulos \cite{revesz/sarantopoulos:2004}.

\medskip Let us settle now our general framework. A \emph{kernel} in the sense
of Fuglede is a lower semicontinuous function $k:X\times
X\rightarrow \RR \cup \{+ \infty\}$ \cite[p.~149]{fuglede:1960}. In
this paper we will sometimes need that the kernel is
\emph{symmetric}, i.e., $k(x,y)=k(y,x)$. This is for example
essential when defining potential and Chebyshev constant, otherwise
there would be a left- and right-potential and the like.

\medskip Another assumption, however a bit of technical flavour, is
the \emph{positivity} of the kernel. This we need, because we would
like to avoid technicalities when integrating not necessarily
positive functions. This assumption is nevertheless not very
restrictive. Since we usually consider  compact sets of $X\times X$,
where by lower semicontinuity $k$ is necessarily bounded from below,
we can assume that $k\ge 0$. Indeed, as we will see, energy,
$n^{\text{th}}$ diameter and $n^{\text{th}}$ Chebyshev constant are
linear in constants added to $k$.

\medskip Denote the set of compactly supported Radon measures on
$X$ by $\fM(X)$, that is
\begin{eqnarray*}
\fM(X):=\{\mu :&\mu& \mbox{ is a regular Borel measure on } X,
\\&\mu& \mbox{ has compact support},\, \|\mu\|<+\infty \}.
\end{eqnarray*}
Further, let $\fMe(X)$ be the set of positive unit measures from
$\fM(X)$,
\begin{equation*}
\fMe(X):= \{\mu \in \fM(X) : \mu\ge 0,\, \mu(X) =1\}.
\end{equation*}
 We say that $\mu\in \fMe(X)$ is supported on $H$ if
$\supp \mu$, which is a compact subset of $X$, is in $H$. The set of
(probability) measures supported on $H$ are denoted by $\fM(H)$
($\fMe(H)$).

\medskip
Before recalling the relevant potential theoretic notions from
\cite{fuglede:1960} (see also \cite{ohtsuka:1961}), let us spend a
few words on integrals (see \cite[Ch.~III-IV.]{bourbaki:1965}). Let
$\mu$ be a positive Radon measure on $X$. Then the integral of a
compactly supported continuous function with respect to $\mu$ is the
usual integral. The upper integral of a positive l.s.c.~function $f$
is defined as
\begin{equation*}
\int\limits_X f\dd\mu:=\sup_{\mbox{\scriptsize$0\leq h\leq f$}\atop
\mbox{\scriptsize$h\in C_c(X)$}} \int\limits_X h \dd\mu.
\end{equation*}
This definition works well, because by standard arguments (see,
e.g., \cite[Ch.~IV., Lemma 1]{bourbaki:1965}) one has
\begin{equation*}
k(x,y)=\sup_{\mbox{\scriptsize$0\leq h\leq k$}\atop
\mbox{\scriptsize$h\in C_c(X\times X)$}} h(x,y),
\end{equation*}
where, because of the symmetry assumption, it suffices to take only
symmetric functions $h$ in the supremum.

What should be here noted, is that this notion of integral has all
useful properties that we are used to in case of Lebesgue integrals
(note also the necessity of the positivity assumptions).

The usual topology on $\fM$ is the so-called \emph{vague topology}
which is a locally convex topology defined by the family
$\{\mu\mapsto \int_X f\dd \mu\::\: f\in C_c(X)\}$ of seminorms. We
will only encounter this topology in connection with families
$\mathcal{M}$ of measures supported on subsets of the same compact
set $K\subset X$. In this case, the weak$^*$-topology (determined by
$C(K)$) and the vague topology coincide on $\mathcal{M}$, Fuglede
\cite{fuglede:1960}.

\medskip For a potential theoretic kernel $k:X\times X\to \RR_+\cup\{0\}$ Fuglede \cite{fuglede:1960} and Ohtsuka \cite{ohtsuka:1961} define the \emph{potential} and the
\emph{energy} of a measure $\mu$
\begin{equation*}
   U^\mu(x) := \int\limits_{X}  k(x,y)\, \dd\mu(y) ~, \qquad
  W(\mu) := \iint\limits_{X\times X}  k(x,y)\, \dd\mu(y)\dd\mu(x).
\end{equation*}
The integrals exist in the above sense, although may attain
$+\infty$ as well.

For a given set $H\subset X$ its \emph{Wiener energy} is
\begin{equation}\label{eq:Energydef}
w(H):= \inf_{\mu\in \fMe(H)} W(\mu),
\end{equation}
see \cite[(2) on p.~153]{fuglede:1960}.

One also encounters the  quantities (see \cite[p.
153]{fuglede:1960})
\begin{equation*}
U(\mu):=\sup_{x\in X} U^\mu(x),\qquad V(\mu):=\sup_{x\in\supp \mu}
U^\mu(x).
\end{equation*}
 Accordingly one defines the following energy functions
\begin{equation*}
u(H):=\inf_{\mu\in\fMe(H)} U(\mu),\qquad v(H):=\inf_{\mu\in\fMe(H)}
V(\mu).
\end{equation*}
In general, one has the relation
\begin{equation*}
 w\leq v\leq u\leq +\infty,
\end{equation*}
where in all places strict inequality may occur. Nevertheless, under
our assumptions we have the equality of the energies $v$ and $w$,
being generally different, see \cite[p.~159]{fuglede:1960}. More
importantly, our set of conditions suffices to have a general
version of Frostman's equilibrium theorem, see Theorem
\ref{thm:frostmantetele}.

In fact, at a certain point (in \S\ref{sec:IM}), we will also assume
Frostman's maximum principle, which will trivially guarantee even
$u=v$, that is, the equivalence of all three energies treated by
Fuglede.

\begin{definition}
The kernel $k$ satisfies the \emph{maximum principle}, if for every
measure $\mu\in \fMe$
\begin{equation*}
  U(\mu)=V(\mu).
 \end{equation*}
\end{definition}
As our examples show in \S\ref{sec:sum}, this is essential also for
the equivalence of the Chebyshev constant and the transfinite
diameter. Carleson \cite[Ch.~III.]{carleson:1967} gives a class of
examples satisfying the maximum principle: Let $\Phi(r)$, $r=|x|$,
$x\in \RR^d$ be the fundamental solution of the Laplace equation,
i.e., $\Phi(|x-y|)$ the Newtonian potential on $\RR^d$. For a
positive, continuous, increasing, convex function $H$ assume also
that
\begin{equation*}
\int\limits_0^1 H(\Phi(r))r ^{d-2}\dd r<+\infty.
\end{equation*}
Then $H\circ \Phi$ satisfies the maximum principle; see
\cite[Ch.~III.]{carleson:1967} and also Fuglede \cite{fuglede:1960}
for further examples. \vskip1ex

Let us now turn to the systematic treatment of the Chebyshev
constant and the transfinite diameter. We call a function $g : X
\rightarrow \RR$ \emph{log-polynomial}, if there exist
$w_1,\ldots,w_n \in X$ such that $g(x)=\sum_{j=1}^n k(x,w_j)$ for
all $x\in X$. Accordingly, we will call the $w_j$s and $n$ the zeros
and the degree of $g(x)$, respectively. Obviously the sum of two
log-polynomials is a log-polynomial again. The terminology here is
motivated by the case of the  logarithmic kernel
\begin{equation*}
k(x,y)=-\log |x-y|,
\end{equation*}
 where the log-polynomials
correspond to negative logarithms of algebraic polynomials.

Log-polynomials give access to the definition of \emph{transfinite
diameter}  and the \emph{Chebyshev constant}, see Carleson
\cite{carleson:1967}, Choquet \cite{choquet:1958}, Fekete
\cite{fekete:1923}, Ohtsuka \cite{ohtsuka:1967} and P\'olya, Szeg\H
o \cite{polya/szego:1931}. First we start with the ``degree $n$''
versions, whose convergence will be proved later.

\begin{definition}
Let $H\subset X$ be fixed. We define the $n^{\text{th}}$ diameter of
$H$ as
\begin{equation}\label{eq:Ddef}
 D_n(H) := \inf_{w_1,\ldots,w_n \in H}
    {\frac{1}{(n-1)\,n}}\bigg(
    \sum_{1\le j \neq l \le n} k(w_j,w_l)
   \bigg);
\end{equation}
or, if the kernel is symmetric
\begin{equation*}
 D_n(H)   =  \inf_{w_1,\ldots,w_n \in H} {\frac{2}{(n-1)\,n}}\bigg(
    \sum_{1\le i <j \le n} k(w_i,w_j) \bigg).
\end{equation*}
\end{definition}
If $H$ is compact, then due to the fact that $k$ is l.s.c., $D_n(H)$
is attained for some points $w_1,\dots,w_n\in H$, which are then
called \emph{$n$-Fekete points}. We will also use the term
\emph{approximate $n$-Fekete points} with the obvious meaning. Note
also that for a finite set $H$, $\#H=m$ and $n>m$, there is always a
point from the diagonal $\Delta=\{(x,x):x\in H\}$ in the definition
of $D_n(H)$. This possibility is completely excluded by Choquet in
\cite{choquet:1958}, thus allowing only infinite sets.
\begin{definition}
For an arbitrary $H \subset X$ the $n^{\text{th}}$ Chebyshev
constant of $H$ is defined as
\begin{equation*}
M_n(H) := \sup_{w_1,\ldots,w_n \in H} \inf_{x\in H}
    {\frac{1}{n}}\bigg(
   \sum_{k=1}^{n} k(x,w_k)
   \bigg)
\end{equation*}
\end{definition}

We are going to show that both $n^{\text{th}}$ diameters and
$n^{\text{th}}$ Chebyshev constants converge from below to some
number (or $+\infty$), which are respectively called the
\emph{transfinite diameter} $D(H)$ and the \emph{Chebyshev constant}
$M(H)$. The aim of this paper is to relate these quantities as well
as the Wiener energy of a set.


\section{Chebyshev constant and transfinite diameter}\label{sec:Cheby}

We define the Chebyshev constant and the transfinite diameter of a
set $H\subset X$ and proceed analogously to the classical case. It
turns out, though not very surprisingly, that in general the
equality of these two quantities does not hold.

First, we prove the convergence of  $n^{\text{th}}$ diameters and
 $n^{\text{th}}$ Chebyshev constants. This is for both cases
classical, we give the proof only for the sake of completeness, see,
e.g., Carleson \cite{carleson:1967}, Choquet \cite{choquet:1958},
Fekete \cite{fekete:1923}, Ohtsuka \cite{ohtsuka:1967} and P\'olya,
Szeg\H o \cite{polya/szego:1931}.

\begin{proposition}
\label{prop:deltanmonoton} The sequence of $n^{\text{th}}$ diameters
is monotonically increasing.
\end{proposition}
\begin{pf}
Choose $x_1,\ldots,x_n\in H$ arbitrarily. If we leave out any index
$i=1,2,\ldots,n$, then for the remaining $n-1$ points we obtain by
the definition of $D_{n-1}(H)$ that
\begin{equation*}
\frac{1}{(n-1)(n-2)}\sum_{\mbox{\scriptsize$1\le j\neq l \le n$}
\atop \mbox{\scriptsize$j\ne i, l\ne i$}}  k(x_j,x_l)  \ge
D_{n-1}(H).
\end{equation*}
After summing up for $i=1,2,\ldots,n$ this yields
\begin{equation*}
\frac{1}{n-1}\sum_{1\leq j\neq l\leq  n} k(x_j,x_l)  \ge n\cdot
D_{n-1}(H),
\end{equation*}
for each term $k(x_j,x_l)$ occurs exactly $n-2$ times. Now taking
the infimum for all possible $x_1,\ldots,x_n \in H$, we obtain
$n\cdot D_n(H) \ge n\cdot D_{n-1}(H)$, hence the assertion.
\qed\end{pf}

\noindent The limit $D(H):=\lim_{n\rightarrow\infty}D_n(H)$ is the
\emph{transfinite diameter} of $H$.

\medskip

\noindent Similarly, the $n^{\text{th}}$ Chebyshev constants
converge, too.
\begin{proposition} For any $H\subset X$, the Chebyshev constants $M_n(H)$
converge in the extended sense.
\end{proposition}
\begin{pf}
The sum of two log-polynomials, $p(z)=\sum_{i=1}^n k(z,x_i)$ with
degree $n$ and $q(z)=\sum_{j=1}^m k(z,y_j)$  with degree $m$, is
also a log-polynomial with degree $n+m$. Therefore
\begin{equation}\label{eq:qm}
{(n+m)}M_{n+m}\ge n M_n + m M_m
\end{equation}
 for all $n,m$ follows at once.
Should $M_n(H)$ be infinity for some $n$, then all succeeding terms
$M_{n'}(H)$, $n'\geq n$ are infinity as well, hence the convergence
is obvious. We assume now that $M_n(H)$ is a finite sequence. At
this point, for the sake of completeness, we can repeat the
classical argument of Fekete \cite{fekete:1923}.

Namely, let $m,n$ be fixed integers. Then there exist $l=l(n,m)$ and
$r=r(n,m),\ 0\le r<m$ nonnegative integers such that $n=l\cdot m +
r$. Iterating the previous inequality \eqref{eq:qm} we get
\begin{equation*}
n\cdot M_n \ge l \Big(m M_m \Big) + r M_r\ = n M_m + r(M_r-M_m).
\end{equation*}
Fixing now the value of $m$, the possible values of $r$ remain
bounded by $m$, and the finitely many values of $M_r-M_m$'s are
finite, too. Hence dividing both sides by $n$, and taking
$\liminf_{n\rightarrow\infty}$, we are led to
\begin{equation*}
\liminf_{n\rightarrow\infty} M_n \ge \liminf_{n\rightarrow\infty}
\bigg( M_m  +  {\frac{r}{n}} \Big( M_r -M_m\Big) \bigg) = M_m \ .
\end{equation*}
This holds for any fixed $m\in\NN$, so taking
$\limsup_{m\rightarrow\infty}$ on the right hand side we obtain
\begin{equation*}
\liminf_{n\rightarrow\infty} M_n \ge \limsup_{m\rightarrow\infty}
M_m ,
\end{equation*}
that is, the limit exists. \qed\end{pf} \noindent
$M(H):=\lim_{n\rightarrow\infty}M_n(H)$ is called the
\emph{Chebyshev constant} of $H$.
\medskip

\noindent In the following, we investigate the connection between
the Chebyshev constant $M(H)$ and the transfinite diameter $D(H)$.
\begin{theorem}\label{thm:Mesdelta}Let $k$ be a positive, symmetric kernel.
 For any $n\in\NN$ and $H\subset X$ we have
$D_n(H)\le M_n(H)$,  thus also $D(H)\le M(H)$.
\end{theorem}
\begin{pf}If $M_n(H)=+\infty$, then the assertion is trivial. So
assume $M_n(H)<+\infty$. By the quasi-monotonicity (see
\eqref{eq:qm}) we have that for all $m\leq n$ also $M_m(H)$ is
finite. We use this fact to recursively find  $w_1,\dots w_n\in H$
such that $k(w_i,w_j)<+\infty$ for all $i<j\leq n$. At the end we
arrive at $\sum_{1\le i<j \le n}k(w_i,w_j)<+\infty$, hence
$D_n(H)<+\infty$. This was our first aim to show, in the following
this choice of the points $w_1,\dots, w_n$ will not play any role.
Instead, for an arbitrarily fixed $\varepsilon>0$, we take, as we
may, an ``approximate $n$-Fekete point system'' $w_1, \dots, w_n$
with
\begin{equation}\label{eq:DnHpe}
\frac{1}{(n-1)n} \sum_{1\le i\neq j \le n}k(w_i,w_j) <
D_n+\varepsilon.
\end{equation}
For any $x\in H$ the points $x,w_1,\ldots,w_n$ form a point system
of $n+1$ points, so by the definition of $D_{n+1}$ we have
\begin{equation*}
  2\sum_{i=1}^{n} k(x,w_i)+ \sum_{1\le i\neq j\le n} k(w_i,w_j) \ge
{n(n+1)} D_{n+1}\ge {{n(n+1)}} D_{n},
\end{equation*}
using also the monotonicity of the sequence $D_n$. This together
with \eqref{eq:DnHpe} lead to
\begin{equation*}
   p_n(x):=\sum_{i=1}^n k(x,w_i) \ge  {\frac{n(n+1)}{2}} D_{n} -
               {\frac{n(n-1)}{2}}\Big( D_{n}+\varepsilon\Big).
\end{equation*}
Taking infimum of the left hand side for $x\in H$ we obtain
\begin{equation*}
   \inf_{x\in H}p_n(x)    \ge n D_n - \frac{n(n-1)\varepsilon}{2}.
\end{equation*}
By the very definition of the $n^{\text{th}}$ Chebyshev constant,
$n\cdot M_n\ge \inf_{x\in H} p_n(x)$ holds, hence $M_n \ge
D_n-(n-1)\varepsilon/2$ follows. As this holds for all
$\varepsilon>0$, we conclude $M_n \ge D_n$. \qed\end{pf}

Later we will show that, unlike the classical case of $\CC$, the
strict inequality $D<M$ is well possible.


\section{Transfinite diameter and energy}\label{sec:DI}
We study the connection between the energy $w$ and the transfinite
diameter $D$. Without assuming the maximum principle we can prove
the equivalence of these two quantities for compact sets. This
result can actually be  found in a note of Choquet
\cite{choquet:1958}. There is however a slight difference to the
definitions of Choquet in \cite{choquet:1958}. There the diagonal
 was completely excluded from the definition of $D$, that is the
infimum in \eqref{eq:Ddef} is taken over $w_i\neq w_j$, $i\neq j$
and not for systems of arbitrary $w_j$'s . This means, among others,
that in \cite{choquet:1958} the transfinite diameter is only defined
for infinite sets. The other assumption of Choquet is that the
kernel is infinite on the diagonal. This is completely the contrary
to what we assume in Theorem \ref{thm:Dwbdk}. Indeed, with our
definitions of the transfinite diameter one can even prove equality
for arbitrary sets if the kernel is finite-valued.

\begin{theorem}\label{thm:DHlewH}
Let $k$ be an arbitrary kernel and $H \subset X$ be any set. Then
$D(H) \le w(H)$.
\end{theorem}
\begin{pf}
Let $\mu\in\fMe(H)$ be arbitrary, and define
$\nu:=\bigotimes_{j=1}^n \mu$ the product measure on the product
space $X^n$. We can assume that the kernel is positive because
$\supp\mu$, and hence $\supp\nu$, is compact  so we can add a
constant to $k$ such that it will be positive on these supports.
Consider the following lower semicontinuous functions $g$ and $h$ on
$X^n$
\begin{eqnarray*}
 g:(x_1,\ldots,x_n)&\mapsto&\displaystyle
D_n(H)\Bigl(:=\inf_{(w_1,\ldots,w_n)\in X^n} {\textstyle \frac
{1}{n(n-1)}} \sum_{1\le i\ne j\le
n}k(w_i,w_j)\Bigr)\\
h:(x_1,\ldots,x_n)&\mapsto&\displaystyle {\textstyle \frac
{1}{n(n-1)}} \sum_{1\le i\ne j\le n} k(x_i,x_j).
\end{eqnarray*}
Since $0\leq g\leq h$, by the definition of the upper integral the
following holds true
\begin{eqnarray*}
D_n(H)& \le& \int\limits_{X^n} \frac {1}{n(n-1)} \sum_{1\le i\ne
j\le n}
k(x_i,x_j) \dd\nu(x_1,\ldots, x_n)  \\
& =& \frac {1}{n(n-1)} \sum_{1\le i\ne j\le n} \int\limits_{H^2}
k(x_i,x_j) \dd\mu(x_i)\dd\mu(x_j)= W(\mu).
\end{eqnarray*}
Taking infimum in $\mu$ yields $D_n(H)\le w(H)$, hence also $D(H)\le
w(H)$. \qed\end{pf}

To establish the converse inequality we need a compactness
assumption. With the slightly different terminology, Choquet proves
the following for kernels being $+\infty$ on the diagonal $\Delta$.
The arguments there are very similar, except that the diagonal
doesn't have to be taken care of in \cite{choquet:1958}. We give a
detailed proof.
\begin{proposition}[Choquet \cite{choquet:1958}]
\label{prop:deltaesi} For an arbitrary kernel function $k$ the
inequality
 $D(K) \ge w(K)$ holds for all $K\subseteq X$ compact sets.
\end{proposition}
\begin{pf}
First of all the l.s.c.~function $k$ attains its infimum on the
compact set $K\times K$. So by shifting $k$ up we can assume that it
is positive, and the validity of the desired inequality is not
influenced by this.

If $D(K)=+\infty$, then by Theorem \ref{thm:DHlewH} we have
$w(K)=+\infty$, thus the assertion follows. Assume therefore
$D(K)<+\infty$, and let $n\in\NN$, $\varepsilon>0$ be fixed. Let us
choose a Fekete point system $w_1,\ldots,w_n$ from $K$. Put
$\mu:=\mu_n:=1/n\, \sum_{i=1}^n \delta_{w_i}$ where $\delta_{w_i}$
are the Dirac measures at the points $w_i$, $i=1,\dots,n$. For a
continuous function $0\leq h\leq k$  with compact support, we have
\begin{eqnarray*}
\iint\limits_{K\times K} h \dd\mu \dd\mu
&=&\frac{1}{n^2}\sum_{i,j=1}^n
h(w_i,w_j)\\
&=&\frac{1}{n^2}\sum_{i=1}^n h(w_i,w_i)+\frac{1}{n^2}\sum_{i,j=1\atop i\neq j}^n h(w_i,w_j)\\
&\leq &\frac{1}{n^2}\sum_{i=1}^n h(w_i,w_i)+\frac{1}{n^2}\sum_{i,j=1\atop i\neq j}^n k(w_i,w_j)\\
&\leq &\frac{\|h\|}{n}+\frac{1}{n^2}\sum_{i,j=1\atop i\neq j}^n k(w_i,w_j)\\
&\leq&\frac{\|h\|}{n}+\frac{n-1}{n}D_n(K)\leq \frac{\|h\|}{n}+D(K)
\end{eqnarray*}
using, in the last step, also the monotonicity of the sequence $D_n$
(Proposition \ref{prop:deltanmonoton}). In fact, we obtain for $n\ge
N= N(\|h\|,\varepsilon)$ the inequality
\begin{equation}\label{eq:IIkmmm}
\iint\limits_{K\times K} h \dd\mu \dd\mu \le D +  \varepsilon.
\end{equation}
It is known, essentially by the Banach-Alaoglu Theorem,  that for a
compact set $K$ the measures of $\fMe(K)$ form a weak$^*$-compact
subset of $\fM$, hence there is a cluster point $\nu\in\fMe(K)$ of
the set $\mathcal{M}_N:=\{\mu_n~:~ n\geq N\}\subset \fMe(K)$. Let
$\{\nu_\alpha\}_{\alpha\in I}\subseteq \mathcal{M}_N$ be a net
converging to $\nu$. Recall that $\nu_\alpha\motimes \nu_\alpha$
weak$^*$-converges to $\nu\motimes\nu$. We give the proof. For a
function $g\in C(K\times K)$, $g(x,y)=g_1(x)\cdot g_2(y)$ it is
obvious that
\begin{equation}\label{eq:wsconv}
\iint\limits_{K\times K} g\dd \nu_\alpha\dd\nu_\alpha\rightarrow
\iint\limits_{K\times K} g\dd \nu\dd\nu.
\end{equation}
The set $\mathcal{A}$ of such product-decomposable functions
$g(x,y)=g_1(x)g_2(y)$ is a subalgebra of $C(K\times K)$, which also
separates $X\times X$, since it is already coordinatewise
separating. By the Stone--Weierstra\ss\ theorem $\mathcal{A}$ is
dense in $C(K\times K)$. From this, using also that the family
$\mathcal{M}_N$ of measures is norm-bounded, we immediately get the
weak$^*$-convergence \eqref{eq:wsconv}. All these imply
\begin{equation*}
\iint\limits_{K\times K}h\dd \nu\dd\nu \leq D(K)+\varepsilon,
\end{equation*}
thus
\begin{equation*}
w(K)\le W(\nu) := \iint\limits_{K\times K} k \dd\nu \dd\nu =
\sup_{\mbox{\scriptsize$0\leq h\leq k$}\atop\mbox{\scriptsize$ h\in
C_c(X\times X)$}}\:\: \iint\limits_{K\times K} h \dd\nu \dd\nu \le
D(K) +  \varepsilon,
\end{equation*}
for all $\varepsilon>0$. This shows $w(K)\le D(K)$. \qed\end{pf}

\begin{corollary}[Choquet \cite{choquet:1958}]\label{cor:DKeqwK} For
arbitrary kernel $k$ and compact set $K\subset X$, the equality
$D(K)=w(K)$ holds.
\end{corollary}
\begin{pf}By compactness we can shift $k$ up and therefore assume
it is positive. Then we apply Theorem \ref{thm:DHlewH} and
Proposition \ref{prop:deltaesi}. \qed\end{pf}
 The assumptions of Choquet \cite{choquet:1958} are the compactness of the set plus the property that the kernel is $+\infty$ on the diagonal (besides it is continuous in
the extended sense). This ensures, loosely speaking, that for a set
$K$ of finite energy
 an energy minimising measure $\mu$ (i.e., for which $W(\mu)=w(K)$) is necessarily non-atomic,
 moreover $\mu\motimes \mu$ is not concentrated on the diagonal.
Therefore to show equality of $w$ with $D$, one has to exclude the
diagonal completely from the definition of the transfinite diameter.

 We however allow a larger set
of choices for the point system in the definition of $D$. Indeed, we
allow Fekete points to coincide, and this also makes it possible to
define the transfinite diameter of finite sets. With this setup the
inequality $D\leq w$ is only simpler than in the case handled by
Choquet. Whereas, however surprisingly, the equality $D(K)=w(K)$ is
still true for compact sets $K$ but without the assumption on the
diagonal values of the kernel.

 We will see in \S\ref{sec:sum} Example
\ref{exa:id} that even assuming the maximum principle but lacking
the compactness allows the strict inequality $D<w$. This phenomena
however may exist only in case of unbounded kernels, as we will see
below. In fact, we show that if the kernel is finite on the
diagonal, then $D=w$ holds for arbitrary sets. For this purpose, we
need the following technical lemma, which shows certain inner
regularity properties of $D$ and is also interesting in itself.

\begin{lemma}\label{lem:DHDKDWbdk}
Assume that the kernel $k$ is positive and finite on the diagonal,
i.e., $ k(x,x) <+\infty$ for all $x\in X$. Then for an arbitrary
$H\subset X$ we have
\begin{equation}\label{eq:DHDKDW}
D(H)=\inf_{\mbox{\scriptsize $K\subset H$}\atop \text{\rm $K$
compact}} D(K) = \inf_{\mbox{\scriptsize $W \subset H$}\atop
\mbox{\scriptsize$\# W<\infty$}} D(W).
\end{equation}
\end{lemma}
\begin{pf} The inequality $\inf D(K)\le\inf
D(W)$ is clear. For $H\supseteq K$ the inequality $D(H)\leq D(K)$ is
obvious, so we can assume $D(H)<+\infty$. For $\varepsilon>0$ let
$W=\{w_1,\ldots,w_n\}$ be an approximate $n$-Fekete point set of $H$
satisfying \eqref{eq:DnHpe}. Then
\begin{equation*}
D(W)=\lim_{m\to\infty} D_{mn}(W)\le \lim_{m\to\infty} \frac
{1}{mn(mn-1)} \sum_{1\le i'\ne j'\le mn} k(w_{i'},w_{j'}),
\end{equation*}
where
\begin{equation*}
 w_{i'}:= \left\{\begin{array}{l} \dots \\ w_i \qquad i'=i+rn,
 \quad r=0,\dots,m-1 \\ \dots
\end{array}\right.
\end{equation*}
Set $C:=\max \{k(x,x):\:x\in W\}$. So we find
\begin{eqnarray*}
  D(W) &\le & \lim_{m\to \infty} \Bigl\{{\textstyle \frac{m^2}{m n(m n-1)}}
  \sum_{1\le i\ne j\le n}\hskip-1em k(w_i,w_j) + {\textstyle\frac{m-1}{m n(m
  n-1)}}
  \sum_{1\le i \le n}\hskip-0.5em k(w_i,w_i)\Bigr\}\\
   &\le & \sum_{1\le i\ne j\le n} k(w_i,w_j) \lim_{m\to\infty}
   {\textstyle \frac{m^2}{m n(m n-1)}} + C n \lim_{m\to\infty} {\textstyle \frac{m-1}{m n(m
   n-1)}}\\
   &=& {\textstyle \frac{1}{n^2}} \sum_{1\le i\ne j\le n} k(w_i,w_j) \le
   {\textstyle \frac{n-1}{n}}\left(D_n(H)+\varepsilon \right) \le D(H)+\varepsilon.
\end{eqnarray*}
This being true for all $\varepsilon>0$, taking infimum we finally
obtain
\begin{equation*}\inf_{\mbox{\scriptsize $W \subset H$}\atop
\mbox{\scriptsize$\# W<\infty$}} D(W)\le D(H).\end{equation*}
\qed\end{pf}

Clearly, if $k(x,x)=+\infty$ for all $x\in W$ with a finite set
$\#W=n$, then for all $m>n$ we have $D_m(W)=+\infty$. Thus in
particular for kernels with $k:\Delta\to\{+\infty\}$, the above can
not hold in general, at least as regards the last part with finite
subsets.

Now, completely contrary to Choquet \cite{choquet:1958} we assume
that the kernel is finite on the diagonal and prove $D=w$ for any
set. Hence an example of $D<w$ (see \S\ref{sec:sum} Example
\ref{exa:id})  must assume $k(x,x)=+\infty$ at least for some point
$x$.

\begin{theorem}\label{thm:Dwbdk}
Assume that the kernel $k$ is positive and is finite on the
diagonal, that is $k(x,x) <+\infty$ for all $x\in X$. Then for
arbitrary sets $H\subset X$, the equality $D(H) = w(H)$ holds.
\end{theorem}
\begin{pf}
 By Theorem \ref{thm:DHlewH} we have $D(H)\leq w(H)$. Hence there is nothing to
prove, if $D(H)=+\infty$. Assume $D(H)<+\infty$, and let
$\varepsilon>0$ be arbitrary. By Lemma \ref{lem:DHDKDWbdk} we have
for some $n\in \NN$ a finite set $W=\{w_1,w_2\ldots,w_n\}$ with
$D(H)+\varepsilon \ge D(W)$. In view of Proposition
\ref{prop:deltaesi} we have $D(W)\ge w(W)$, and by monotonicity also
$w(W)\ge w(H)$. It follows that $D(H)+\varepsilon \ge w(H)$ for all
$\varepsilon>0$, hence also the ``$\ge$'' part of the assertion
follows. \qed\end{pf}


\section{Energy and Chebyshev constant}\label{sec:IM}
To investigate the relationship between the energy and the Chebyshev
constant the following general version of Frostman's Equilibrium
Theorem \cite[Theorem 2.4]{fuglede:1960} is fundamental for us.
\begin{theorem}[Fuglede]
\label{thm:frostmantetele} Let $k$ be a positive, symmetric kernel
and $K\subset X$ be a compact set such that $w(K)<+\infty$. Every
$\mu$ which has minimal energy ($\mu \in \fMe(K), W(\mu)=w(K)$)
satisfy the following properties
\begin{eqnarray*}
U^\mu(x) & \ge& w(K)\quad\mbox{ for nearly every\footnotemark}\ x \in K,\\
U^\mu(x) & \le& w(K) \quad\mbox{ for every } x \in \supp \mu, \\
 U^\mu(x) & =& w(K)\quad \mbox{ for } \mu\mbox{-almost every }
x \in X.
\end{eqnarray*}
Moreover, if the kernel is \emph{continuous}, then
\begin{eqnarray*}
U^\mu(x) & \ge w(K)\quad\mbox{ for \emph{every} }x \in
K.\rule{3.4em}{0pt}
\end{eqnarray*}\footnotetext{The set $A$ of exceptional points is small in the sense $w(A)=+\infty$.}
\end{theorem}

\begin{theorem}
\label{thm:iesm} Let $H\subset X$ be arbitrary. Assume that the
kernel $k$ is positive, symmetric and satisfies the maximum
principle. Then we have $M_n(H)\le w(H)$ for all $n\in\NN$, whence
also $M(H)\le w(H)$ holds true.
\end{theorem}
\begin{pf}
Let $n\in\NN$ be arbitrary. First let $K$ be any compact set. We can
assume $w(K)<+\infty$, since otherwise the inequality holds
irrespective of the value of $M_n(K)$. Consider now an
energy-minimising measure $\nu_K$ of $K$, whose  existence is
assured by the lower semicontinuity of $\mu \mapsto \iint
k\dd\mu\dd\mu$ and the compactness of $\fMe(K)$, see \cite[Theorem
2.3]{fuglede:1960}.

By the Frostman-Fuglede theorem (Theorem \ref{thm:frostmantetele})
we have $ U^{\nu_K}(x) \le w(K)$ for all $x \in \supp \nu_K$, so
$V(\nu_K)\leq w(K)$,
   and by the maximum principle even
   \begin{equation*}
U^{\nu_K}(x) \le w(K) \quad \mbox{ for all } x \in
   X.
\end{equation*}
Then for all $w_1,\ldots,w_n \in K$
\begin{equation*}
 \inf_{x\in K} \frac{1}{n}\sum_{j=1}^n k(x,w_j) \le
\int\limits_X \frac{1}{n}\sum_{j=1}^n k(x,w_j) \dd\nu_K(x) \le w(K)\
.
\end{equation*}
Taking supremum for $w_1,\ldots,w_n \in K$, we obtain
\begin{equation*}
  \sup_{w_1,\ldots,w_n \in K} \inf_{x\in K} \frac{1}{n}\sum_{j=1}^n k(x,w_j)
\le w(K).
\end{equation*}
So $M_n(K)\le w(K)$ for all $n\in\NN$.

Next let $H\subset X$ be arbitrary. In view of the last form of
\eqref{eq:Energydef}, for all $\varepsilon>0$ there exists a measure
$\mu\in\fMe(H)$, compactly supported in $H$, with $w(\mu)\le
w(H)+\varepsilon$. Let $W=\{w_1,\dots, w_n\}\subset H$ be arbitrary
and define $p_W(x):=\frac{1}{n}\sum_i k(x,w_i)$.

Consider the compact set $K:=W\cup\supp \mu\subset H$. By definition
of the energy, $\supp\mu\subset K$ implies $w(K)\le w(\mu)$, hence
$w(K)\le w(H)+\varepsilon$. Combining this with the above, we come
to $M_n(K)\le w(H)+\varepsilon$. Since $W\subset K$, by definition
of $M_n(K)$ we also have
\begin{equation}\label{eq:pWMnK}
\inf_{x\in K} p_W(x) \le M_n(K).
\end{equation}
The left hand side does not increase, if we extend the $\inf$ over
the whole of $H$, and the right hand side is already estimated from
above by $w(H)+\varepsilon$. Thus \eqref{eq:pWMnK} leads to
\begin{equation*}
\inf_{x\in H} p_W(x) \le w(H)+\varepsilon.
\end{equation*}
This holds for all possible choices of $W=\{w_1,\dots,w_n\}\subset
H$, hence is true also for the $\sup$ of the left hand side. By
definition of $M_n(H)$ this gives exactly $M_n(H)\le
w(H)+\varepsilon$, which shows even $M_n(H)\le w(H)$. \qed\end{pf}

\begin{remark*}
In \cite{farkas/revesz:2004c} it is proved that $M(H)=q(H)$, where
\begin{equation*}
q(H)=\inf_{\mu\in\fMe(H)}\sup_{x\in H} U^\mu(x).
\end{equation*}
The idea behind is a minimax theorem, see also
\cite{ohtsuka:1965,ohtsuka:1967}.
 Trivially $w(H)\leq
q(H)\leq u(H)$. So the maximum principle implies
$M(H)=w(H)=q(H)=u(H)$.
\end{remark*}


\section{Summary of the Results. Examples}\label{sec:sum}
In this section, we put together the previous results, thus proving
the equality of the three quantities being studied, under the
assumption of the maximum principle for the kernel. Further, via
several instructive examples we investigate the necessity of our
assumptions and the sharpness of the results.

\begin{theorem}\label{thm:sumth} Assume that the kernel $k$ is positive, symmetric and
satisfies the maximum principle. Let $K\subset X$ be any compact
set. Then the transfinite diameter, the Chebyshev constant and the
energy  of $K$ coincide:
\begin{equation*}
D(K)=M(K)=w(K).
\end{equation*}
\end{theorem}
\begin{pf} We presented a cyclic proof above, consisting
of $M\ge D$ (Theorem \ref{thm:Mesdelta}), $D\ge w$ (Proposition
\ref{prop:deltaesi}) and finally  $w\ge M$ (Theorem \ref{thm:iesm}).
\qed\end{pf}

\begin{theorem}\label{thm:sumthbdd} Assume that the kernel $k$ is positive, finite
and satisfies the maximum principle. For an arbitrary subset
$H\subset X$ the transfinite diameter, the Chebyshev constant and
the energy  of $H$ coincide:
\begin{equation*}
D(H)=M(H)=w(H).
\end{equation*}
\end{theorem}
\begin{pf} By finiteness $D=w$, due to Theorem \ref{thm:Dwbdk}. This with $D\leq M$ and $M\leq w$
(Theorems \ref{thm:Mesdelta} and \ref{thm:iesm}) proves the
assertion. \qed\end{pf}
\begin{remark*}
In the above theorem, logically it would suffice to assume that the
kernel be finite only on the diagonal. But if this was the case, the
maximum principle would then immediately imply the finiteness of the
kernel \emph{everywhere}.
\end{remark*}
\bigskip
Let us now discuss how sharp the results of the preceding sections
are. In the first example we show that, if we drop the assumption of
compactness the assertions of Theorem \ref{thm:Mesdelta}, Theorem
\ref{thm:DHlewH} and Theorem \ref{thm:iesm} are in general the
strongest possible.

\begin{example}\label{exa:id}
Let $X=\NN\cup\{0\}$ endowed with discrete topology and the kernel
\begin{equation*}
k(n,m):=\left\{\begin{array}{ccl}
+\infty&\quad\quad&\mbox{if $n=m$},\\
0&&\mbox{if $0\neq n\neq m\neq 0$},\\
1&& \mbox{otherwise}.
\end{array}\right.
\end{equation*}
The kernel is symmetric, l.s.c.~and has the maximum principle. This
latter can be seen by noticing that for a probability measure
$\mu\in \fMe(X)$ the potential is $+\infty$ on the support of $\mu$.
Indeed,  since $X$ is countable, all measures $\mu\in \fMe(X)$ are
necessarily atomic, and if for some point $\ell\in X$ we have
$\mu(\{\ell\})>0$, then by definition $\int_X k(x,y) \dd\mu(y)
=+\infty$.

We calculate the studied quantities of the set $H=X$ (also as in all
the examples below). Since the kernel is positive, $D_n\geq 0$. On
the other hand, choosing $w_1:=1,\ldots,w_n:=n$, all the values
$k(w_i,w_j)$ will be exactly $0$, so it follows that $D_n=0$,
$n=1,2,\ldots$, and hence $D=0$.

The Chebyshev constant can be estimated from below, if we compute
the infimum of a suitably chosen log-polynomial.  Consider the
log-polynomial $p(x)$ with all zeros placed at $0$, that is with
$w_1=\ldots=w_{n}=0$. Then the log-polynomial $p(x)$ is $\sum_j
k(x,w_j)=n\cdot k(x,0)$. If $x\neq 0$, we have $p(x)= n$, which
gives $M\geq 1$. The upper estimate of $M$ is also easy: suppose
that in the system $w_1,\ldots w_n$ there are exactly $m$ points
being equal to $0$ (say the first $m$). Then
\begin{equation*}
p(x)=\left\{\begin{array}{ccl}
+\infty &\quad\quad&x=w_1,\ldots,w_n, \\
n&\quad\quad& x=0,\: x\neq w_1,\ldots,w_n\quad\quad\mbox{(if $m=0$)}\\
m&\quad\quad& x\neq0,\: x\neq w_1,\ldots,w_n
\end{array}\right.
\end{equation*}
This shows for the corresponding log-polynomial $\inf p(x)=m$, so
$M_n\leq 1$, whence $M=1$.

The energy is computed easily. Using the above reasoning on the
maximum principle, we see $W(\mu)=+\infty$ for any $\mu\in \fMe(X)$,
hence $w(X)=+\infty$.

Thus we have an example of
\begin{equation*}
+\infty=w>M>D=0.
\end{equation*}
\end{example}
The above example completes the case of the kernel with maximum
principle. Let us now drop this assumption and look at what can
happen.

\begin{example}
Let $X:=\{-1,0,1\}$ be endowed with the discrete topology. We define
the kernel by
\begin{equation*}
k(x,y) := \left\{
\begin{array}{ll}
2\quad\quad & \mbox{ if } 0\le|x-y|< 2,\\
0\quad\quad & \mbox{ if } 2=|x-y|.
\end{array}
\right.
\end{equation*}
Then $k$ is continuous and bounded on $X\times X$. This, in any
case, implies $D=w$ by Theorem \ref{thm:Dwbdk}. Note that $k$ does
not satisfy the maximum principle. To see this, consider, e.g., the
measure $\mu=\frac{1}{2}\delta_{-1}+\frac{1}{2}\delta_1$. Then for
the potential $U ^\mu$ one has $U ^\mu(1)=U ^\mu(-1)=1$ and $U
^\mu(0)=2$, which shows the failure of the maximum principle.

To estimate the $n^{\text{th}}$ diameter from above, let us consider
the point system $\{w_i\}$ of $n=2m$ points with $m$ points falling
at $-1$ and $m$ points falling at $1$, while no points being placed
at $0$. Then by definition of $D_n:=D_n(X)$ one can write
\begin{equation*}
\frac{n(n-1)}{2}D_n \le 2 \left(m \atop 2 \right)\cdot 2 + m^2\cdot
0 = \frac{n^2}{2} - n.
\end{equation*}
Applying this estimate for all even $n=2m$ as $n\to\infty$, it
follows that
\begin{equation}\label{eq:DDDD}
D=\lim_{n\to\infty} D_n \le 1.
\end{equation}

Next we estimate the Chebyshev constants from below by computing the
infimum of some special log-polynomials. For $p_n(x)= k(x,0)$ one
has $p_n(x)\equiv 2=\inf p_n$. We thus find $M_n \ge 2$ and $M \ge
2$, showing $M>D$, as desired.
\end{example}

\begin{example}\label{exa:1} Let $X:=\NN$ with the discrete topology. Then $X$ is a
locally compact Hausdorff space, and all functions are continuous,
hence l.s.c.~on $X$. Let $k:X\times X\to [0,+\infty]$ be defined as
\begin{equation*}
k(n,m):=\left\{\begin{array}{ll}+\infty &{\rm if} \qquad n=m, \\
2^{-n-m} & \mbox{if} \qquad n\ne m.
\end{array}\right.
\end{equation*}
Clearly $k$ is an admissible kernel function. For the energy we have
again $w(X)=+\infty$, see Example \ref{exa:id}.

On the other hand let $n\in\NN$ be any fixed number, and compute the
$n^{\rm th}$ diameter $D_n(X)$. Clearly if we choose $w_j:=m+j$,
with $m$ a given (large) number to be chosen, then we get
\begin{equation*}
D_n(H)\le {\frac{1}{(n-1)\,n}}\sum_{1\le i \ne j \le n} 2^{-i-j-2m}
\le \frac{2^{-2m}}{(n-1)\,n} \bigg(\sum_{i=1}^\infty
    2^{-i}\bigg)^2\le 2^{-2m}~,
\end{equation*}
hence we find that the $n^{\rm th}$ diameter is $D_n(X)=0$, so
$D(X)=0$, too. For any log-polynomial $p(x)$ we have $\inf
p(x)=\lim_{x\to \infty} p(x)=0$, hence $M(X)=0$.  That is we have
$D(X)=M(X)=0<w(X)=+\infty$.
\end{example}

The example shows how important the diagonal, excluded in the
definition of $D$ but taken into account in $w$, may become for
particular cases. We can even modify the above example to get finite
energy.
\begin{example}
Let $X:=(0,1]$, equipped with the usual topology, and let $x_n=1/n$.
We take now
\begin{equation*}
k(x,y):=\left\{\begin{array}{ll}
+\infty &\mbox{if $x=y$}, \\
2^{-n-m} \qquad \qquad  &\mbox{if  $x=x_n \mbox{ and } \,\; y=x_m \;\,(x_n\ne x_m)$},\\
-\log|x-y|&\mbox{otherwise}
\end{array}\right.
\end{equation*}
Compared to the l.s.c.~logarithmic kernel, this $k$ assumes
different, smaller values at the relatively closed set of points $\{
(x_n,x_m)~:~ n \ne m \}\subset X\times X$ only, hence it is also
l.s.c.~and thus admissible as kernel.

If a measure $\mu\in \fMe(X)$ has any atom, say if for some point
$z\in X$ we have $\mu(\{z\})>0$, then by definition $\int_X k(x,y)
\dd\mu(y) =+\infty$, hence also $w(\mu)=+\infty$. Since for all
$\mu\in \fMe(X)$ with any atomic component $w(\mu)=+\infty$, we find
that for the set $H:=X$ we have
\begin{equation*}
w(H):=\inf_{\mu\in \fMe(H)} w(\mu) = \inf_{\mu\in \fMe(H) \atop {
\text{$\mu$ not atomic}}} w(\mu).
\end{equation*}
But for measures without atoms, the countable set of the points
$x_n$ are just of measure zero, hence the energy equals to the
energy with respect to the logarithmic kernel. Thus we conclude
$w(H)=e^{-{\rm cap}(H)}=e^{-1/4}$, as ${\rm cap}((0,1])=1/4$ is
well-known.

On the other hand if $n\in\NN$ is any fixed number, we can compute
the $n^{\rm th}$ diameter $D_n(H)$ exactly as above in Example
\ref{exa:1}. Hence it is easy to see that $D_n(H)=0$, whence also
$D(H)=0$. Similarly, we find $M(H)=0$, too.

This example shows that even in case $w(H)<+\infty$ we can have
$w(H)>D(H)=M(H)$.
\end{example}


\section{Average distance number and the maximum principle}
In the previous section, we showed the equality of the Chebyshev
constant $M$ and the transfinite diameter $D$, using essentially
elementary inequalities and the only theoretically deeper
ingredient, the assumption of the maximum principle. We have also
seen examples showing that the lack of the maximum principle for the
kernel allows strict inequality between $M$ and $D$. These
observations certify to the relevance of this principle in our
investigations. Indeed, in this section we show the necessity of the
maximum principle in case of  continuous kernels for having
$M(K)=D(K)$ for all compact sets $K$. We need some preparation
first.

\medskip
Recall from the introduction the notion of the average distance (or
rendezvous) number. Actuyally, a more general assertion than there
can be stated, see Stadje \cite{stadje:1981} or
\cite{farkas/revesz:2004c}. For a compact connected set $K$ and a
continuous, symmetric kernel $k$, the average distance number $r(K)$
is the uniquely existing number with the property that for all
probability measures supported in $K$ there is a point $x\in K$ with
\begin{equation*}
U^\mu(x)=\int\limits_Kk(x,y)\dd\mu(y)=r(K).
\end{equation*}
This can be even further generalised by dropping the connectedness,
see Thomassen \cite{thomassen:2000} and \cite{farkas/revesz:2004c}.
Even for not necessarily connected but compact spaces $K$ with
symmetric, continuous kernel $k$ there is a unique number $r(K)$
with the property that whenever a probability measure on $K$ and a
positive $\varepsilon$ are given, there are points $x_1,x_2\in K$
such that
\begin{equation*}
U^\mu(x_1)-\varepsilon\leq r(K)\leq U^\mu(x_2)+\varepsilon.
\end{equation*}
This number is called the \emph{(weak) average distance number}, and
is particularly easy to calculate, when a probability measure with
constant potential is available. Such a measure $\mu$ is called then
an \emph{invariant measure}. In this case the average distance
number $r(K)$ is trivially just the constant value of the potential
$U^\mu$, see Morris, Nicholas \cite{morris/nicholas:1983} or
\cite{farkas/revesz:2004b}.

It was proved in \cite{farkas/revesz:2004b} that one always has
$M(K)=r(K)$, so once we have an invariant measure, then the
Chebyshev constant is again easy to determine.

Also the Wiener energy $w(K)$ has connection to invariant measures,
as shown by the following result, which is a simplified version of a
more general statement from \cite{farkas/revesz:2004b}, see also
Wolf \cite{wolf:1997}.
\begin{theorem}\label{thm:kinvariant} Let $\emptyset\neq K\subset X$ be a
compact set and $k$ be a continuous, symmetric kernel. Then we have
\begin{equation*}
r(K)\geq w (K).
\end{equation*}
 Furthermore, if $r(K)= w(K)$, then there exists an
invariant measure in $\fMe(K)$.
\end{theorem}
As mentioned above, we have $r(K)=M(K)$, so the inequality $r(K)\geq
w(K)$ in the first assertion of the above theorem is also the
consequence of Theorems \ref{thm:Mesdelta} and \ref{thm:Dwbdk}. For
the proof of the second assertion one can use the Frostman-Fuglede
Equilibrium Theorem \ref{thm:frostmantetele} with the obvious
observation that ``nearly every'' in this context means indeed
``every''. Actually any probability measure $\mu\in \fMe(K)$ which
minimises $\nu\mapsto \sup_KU^\nu$ is an invariant measure and its
potential is constant $M(K)$, see
\cite[Thm.~5.2]{farkas/revesz:2004b} (such measures undoubtedly
exist because of compactness of $\fMe(K)$). Henceforth we will
indifferently use the terms energy minimising or invariant for
expressing this property of measures.
\begin{theorem}
Suppose that the kernel $k$ is symmetric and continuous. If
$M(K)=D(K)$ for all compact sets $K\subseteq X$, then the kernel has
the maximum principle.
\end{theorem}
\begin{pf}
Recall from Corollary \ref{cor:DKeqwK} that $D(K)=w(K)$ for all
$K\subseteq X$ compact. So we can use Theorem \ref{thm:kinvariant}
all over in the following arguments. We first prove the assertion in
the case when $X$ is a finite set. The proof is by induction on
$n=\# X$. For $n=1$ the assertion is trivial. Let now $\# X=2$,
$X=\{a,b\}$. Assume without loss of generality that $k(a,a)\leq
k(b,b)$. Then we only have to prove that for $\mu=\delta_a$ the
maximum principle, i.e., the inequality $k(a,b)\leq k(a,a)$ holds.
To see this we calculate $M(X)$ and $D(X)$. We certainly have
$D(X)\leq k(a,a)$. On the other hand for an energy minimising
probability measure $\nu_p:=p \delta_a+(1-p)\delta_b$ on $X$ we know
that its potential is constant over $X$, hence
\begin{eqnarray*}
p k(a,a)+(1-p)k(b,a)&=&pk(a,b)+(1-p)k(b,b)\\
&=&M(X)=D(X)\leq k(a,a).
\end{eqnarray*}
Here if $p=1$, then $k(a,a)=k(a,b)$. If $p<1$, then we can write
\begin{equation*}
(1-p)k(b,a)\leq  (1-p)k(a,a),\quad\mbox{hence}\quad k(b,a)\leq
k(a,a),
\end{equation*}
 so the maximum principle holds.

 Assume now that the assertion is true for all sets with at most $n$
 elements and for all kernels, and let $\#X=n+1$.
 For a probability measure $\mu$ on $X$ we have to prove $\sup_{x\in X}
 U^\mu(x)=\sup_{x\in\supp\mu} U^\mu(x)$. If $\supp \mu=X$, then there is
nothing to prove. Similarly, if there are two distinct points
$x_1\neq x_2$, $x_1,x_2\in X\setminus\supp\mu$, then by the
induction hypothesis we have
\begin{equation*}
\sup_{x\in X\setminus\{x_1\}} U^\mu(x)=\sup_{x\in\supp \mu}U^\mu(x)=
\sup_{x\in X\setminus\{x_2\}} U^\mu(x).
\end{equation*}
So for a probability measure $\mu$ defying the maximum principle we
must have $\#\supp\mu=n$, say $\supp\mu=X\setminus\{x_{n+1}\}$; let
$\mu$ be such a measure. Set $K=\supp\mu$ and let $\mu'$ be an
invariant measure on $K$.  We claim that all such measures $\mu'$
are also violating the maximum principle. If $\mu=\mu'$, we are
done. Assume $\mu\neq \mu'$ and consider the linear combinations
$\mu_t:=t\mu +(1-t)\mu'$. There is a $\tau>1$, for which
$\mu_{\tau}$ is still a probability measure and
$\supp\mu_{\tau}\subsetneq \supp\mu$. By the inductive hypothesis
(as $\#\supp \mu_\tau<n$) we have $U^{\mu_{\tau}}(x_{n+1})\leq
U^{\mu_{\tau}}(a)$ for some $a\in \supp \mu_\tau$. We also know that
$U^\mu(x_{n+1})=U^{\mu_1}(x_{n+1})> U^{\mu_1}(a)$. Hence for the
linear function $\Phi(t):=U^{\mu_t}(x_{n+1})-U^{\mu_t}(a)$ we have
$\Phi(1)>0$ and also  $\Phi(\tau)\leq 0$ ($\tau>1$). This yields
$\Phi(0)>0$, i.e.,  $U^{\mu'}(x_{n+1})=U^{\mu_{0}}(x_{n+1})>
U^{\mu_{0}}(a)=U^{\mu'}(y)$ for all $y\in K$.
 We have therefore shown that all energy minimising
(invariant) measures on $K$ must defy the maximum principle.

 Let now $\nu$  be an
invariant measure on $X$. We have
\begin{eqnarray*}
M(X)&=&U^\nu(y)=\sup_{x\in X} U^\nu(x)=D(X)\\
&\leq& D(K)=\sup_{x\in K} U^{\mu'}(x)=U^{\mu'}(z)< U^{\mu'}(x_{n+1})
\end{eqnarray*}
for all $y\in X$, $z\in K$. Thus we can conclude $U^\nu(y)\leq
U^{\mu'}(y)$ for all $y\in X$ and even ``$<$'' for $y=x_{n+1}$.
Integrating with respect to $\nu$ would yield
\begin{equation*}
\int\limits_X\int\limits_X k\dd\nu\dd\nu
=M(X)<\int\limits_X\int\limits_Xk
\dd\mu'\dd\nu=\int\limits_X\int\limits_Xk \dd\nu\dd\mu'=M(X),
\end{equation*}
hence a contradiction, unless $\nu(\{x_{n+1}\})=0$. If
$\nu(\{x_{n+1}\})=0$ held, then $\nu$ would be an energy minimising
measure on $K$. This is because obviously $\supp\nu\subseteq K$
holds, and the potential of $\nu$ is constant $M(X)$ over $K$, so
\begin{equation*}
M(X)=\int\limits_K\int\limits_K
k\dd\nu\dd\mu'=\int\limits_K\int\limits_K
k\dd\mu'\dd\nu=M(K)\quad\mbox{holds}.
\end{equation*}
 As we saw above, then $\nu$ would not satisfy the maximum principle, a contradiction
again, since the potential of $\nu$ is constant on $X$. The proof of
the case of finite $X$ is complete.

We turn now to the general case of $X$ being a locally compact space
with continuous kernel. Let $\mu$ be a compactly supported
probability measure on $X$ and $y\not\in \supp \mu$. Set $K=\supp
\mu$ and note that both $\fMe(K)\ni \nu \mapsto\sup_K U^\nu$ and
$\nu\mapsto U^\nu (y)$ are continuous  mappings with respect to the
weak$^*$-topology on $\fMe(K)$. If $\sup_K U^\mu<U^\mu(y)$ were
true, we could therefore find, by a standard approximation argument,
see for example \cite[Lemma 3.8]{farkas/revesz:2004c}, a finitely
supported probability measure $\mu'$ on $K$ for which
\begin{equation*}
\sup_{x\in\supp\mu'} U^{\mu'}(x)\leq \sup_{x\in K}
U^{\mu'}(x)<U^{\mu'}(y).
\end{equation*}
 This is nevertheless impossible by the first part of the proof,
 thus  the assertion of the theorem follows.
\qed\end{pf}


\section*{Acknowledgement}
The authors are deeply indebted to Szil\'ard R\'ev\'esz for his
insightful suggestions and for the motivating discussions.


\end{document}